\documentclass[12pt]{article}%
\usepackage{amsmath}
\usepackage{amsfonts}
\usepackage{amssymb}
\usepackage{graphicx}%
\usepackage{mathrsfs}

\usepackage{amssymb}
\usepackage{stmaryrd}
\usepackage{tipa}
\usepackage{amsfonts,amsmath,latexsym}
\usepackage{mathrsfs}  
\usepackage{amsbsy}
\usepackage{indentfirst}
\usepackage{amsmath}
\usepackage{amsfonts}
\usepackage{bbm}
\usepackage{mathrsfs}
\usepackage{dsfont}
\usepackage{graphicx}

\allowdisplaybreaks[4]

\usepackage[square, comma, sort&compress, numbers]{natbib}
\usepackage{amssymb, color}
\usepackage[body={15.5cm,21cm}, top=3cm]{geometry}%
\usepackage{hyperref}
\providecommand{\U}[1]{\protect \rule{.1in}{.1in}}

\newtheorem{theorem}{Theorem}[section]

\newtheorem{definition}[theorem]{{Definition}}

\newtheorem{lemma}[theorem]{Lemma}

\newtheorem{remark}[theorem]{{Remark}}

\begin{document}

\title{CENTRAL LIMIT THEOREM  FOR IRREGULAR DISCRETIZATION SCHEME OF MULTILEVEL MONTE CARLO METHOD}
\author{ Yi Guo
, Yuxi Guo
,  Hanchao Wang\\
Institute for Financial Studies, Shandong University,  Jinan,  250100, China
}
\linespread{1.5}
\definecolor{shadecolor}{RGB}{241, 241, 255}

\maketitle

\begin{abstract}
In this paper, we study the asymptotic error distribution for a two-level irregular discretization scheme of the solution to the stochastic differential equations (SDE for short) driven by a continuous semimartingale and obtain a central limit theorem for the error processes with the rate $\sqrt{n}$. As an application, in the spirit of the result of  Ben Alaya and Kebaier, we get a central limit theorem of the Linderberg-Feller type for the irregular discretization scheme of the multilevel Monte Carlo method.
\end{abstract}

\medskip
\textbf{Key words}:  Stable convergence, irregular discretization scheme, multilevel Monte Carlo scheme, central limit theorem.

\medskip
\noindent \textbf{MSC-classification}:  60F17, 60H35.

\section{Introduction}
The Monte Carlo method is an important technique in many applications. For instance, we are interested in the law of $X_T$, where $\{X_t\}_{0\le t\le T}$ is a given process. However, in most cases, the law of $X_T$ is unknown, and we are expected to compute the $\mathbb{E}[f(X_T)]$ for some test functions. The Monte Carlo method is then applied to approximate the $\mathbb{E}[f(X_T)]$ by the estimator 
\[Z_n=\frac{1}{N}\sum_{j=1}^Nf(X_{T,j})\]
where $N$ is a constant and $\{X_{T,j}\}_{1\le j\le N}$ are N independent copies of $X_T$. Moreover, if $X$ is a diffusion process, we usually use the Euler scheme $X^n$ with time step $T/n$ to approximate $X$.\par
To improve the accuracy of the classic Monte Carlo method, Kebaier first proposed a two-level Monte Carlo Euler method in \cite{kebaier2005statistical}. Giles \cite{giles2008multilevel} extended this method and introduced the multilevel Monte Carlo method. To be more precise, the quantity $\mathbb{E}[f(X_T)]$ is estimated by
\[Q_n=\frac{1}{N_0}\sum_{k=1}^{N_0}f(X_{T,k}^1)+\sum_{\ell=1}^L\frac{1}{N_{\ell}}\sum_{k=1}^{N_{\ell}}(f(X_{T,k}^{m^{\ell}})-f(X_{T,k}^{m^{\ell-1}})),\]
where $L=\log n/\log m$ and for fixed $\ell\in\{0,\cdots,L\}$, $X_{t}^{m^{\ell}}$ is the Euler scheme with time step $m^{-\ell}T$ and $\{X_{T,k}^{m^{\ell}}\}_{1\le k\le N_{\ell}}$ are $N_{\ell}$ independent copies of $X_{t}^{m^{\ell}}$. It is worth noting that the simulation of all $\{X_{t}^{m^{\ell}}\}_{1\le \ell\le L}$ should be based on the same path of Brown motion. Ben Alaya and Kebaier \cite{ben2015central} gave a central limit theorem (see Theorem 4 of \cite{ben2015central}) for the error of the multilevel Monte Carlo method based on a stable convergence theorem of $\frac{nm}{(m-1)T}(X_{T,k}^{m^{\ell}}-X_{T,k}^{m^{\ell-1}})$ which seems to be a generalization of Jacod and Protter \cite{jacod98}. They proved that the convergence rate can be $1/n^{\alpha}$ for all $\alpha\in[1/2,1]$. It turns out that this multilevel method has more applications. For example, Ben Alaya et al. \cite{arxiv} obtained the asymptotic error of the multilevel Euler method for SDEs driven by a pure  L\'evy process, which is an extension of Jacod \cite{jacod04} and Wang \cite{wang2015euler}.
\par
In effect, the authors above are based on the regular Euler scheme, at which point the sampling is equidistant, and we realize that it is more important when the scheme is not regular. Therefore, in this paper, we focus on the irregular discretization scheme, which amounts to say that $X$ is sampled at $\tau_{k}^n$ for $k=1,2,\cdots$ which have the form
\begin{equation}\label{time}
\tau_{k+1}^n=(\tau_{k}^n+\frac{1}{n\theta(\tau_{k}^n)}),
\end{equation}
where $\theta$ is an adapted process and $|1/\theta|$ is bounded by a constant $K$. It can be reduced to a regular scheme if we let $\theta=1/T$. For more details on the irregular scheme, readers can refer to Jacod and Protter \cite{jacod2012discretization}, Part \uppercase\expandafter{\romannumeral6}. \par
Inspired by Ben Alaya and Kebaier \cite{ben2015central}, we will prove a central limit theorem for multilevel irregular discretization schemes. Our main result (see Theorem \ref{main}) extends Theorem 3 in Ben Alaya and Kebaier \cite{ben2015central}, at which point the scheme is regular. The conditions required to converge the error processes are more stringent than that paper. Specifically, we need $1/\theta$ to be almost surely Lebesgue integrable over $[0, T]$ and to approximate the quadratic variation and covariation processes. The approximation tools are given in Lemma \ref{carl}, which can be referred to  Lindberg \cite{lindberg2013error}. As an application, a central limit theorem for the multilevel Monte Carlo method is obtained. \par
The rest of the paper is organized as follows: Section 2 first introduces the general settings and defines some notations. The main theorem is presented in Section 3, and we divide the proof of the theorem into three steps, with the third step being the most important. In Section 4, we derive a central limit theorem for the multilevel Monte Carlo method and an application to call option using our main theorem.

\section{General framework.}
\subsection{Preliminaries}
Let $(X^n)_{n\ge 1}$ be a sequence of random variables defined on some probability space $(\Omega, \mathcal{F},  \mathbb{P})$ and with values in $C[0,T]$. Let $(\widetilde{\Omega},\widetilde{\mathcal{F}}, \widetilde{ \mathbb{P}})$ be an extension of $(\Omega, \mathcal{F},  \mathbb{P})$. Then we say $(X^n)_{n\ge 1}$ converges stably to $X$ if for all bounded continuous function $f:C[0,T]\to \mathbb{R}$ and all bounded measurable random variable $U$, 
\[\mathbb{E}[Uf(X^n)]\to \widetilde{ \mathbb{E}}[Uf(X)].\]
Then we denote stable convergence by $\stackrel{s-\mathcal{L}}{\longrightarrow}$. Weak convergence is denoted by $\Rightarrow$.\par
 In the following, we recall the PUT condition and two theorems, which can be found in \cite{jacod2013} \uppercase\expandafter{\romannumeral6}.6.15, 6.22 and are crucial for the convergence of the stochastic integrals.
\begin{definition}
A sequence of  continuous $\mathbb{R}^d$-valued semimartingales $(X^n)_{n\ge 1}$ is said to have PUT condition if for each $n\ge 1$, the decomposition $X^n= M^n+A^n$ have
\[\sup_n\mathbb{E}[\left\langle M^n ,M^n   \right\rangle_T+\int_0^T|dA_s^n|]<\infty.\]
\end{definition}
\begin{theorem}\label{put}
Let $H^n$ and $H$ be a sequence of adapted,  $c\grave{a}dl\grave{a}g$ $\mathbb{R}^{q\times d}$ processes defined on probability space where $c\grave{a}dl\grave{a}g$ means right continuous with left limits. Assume $(X^n)_{n\ge 1}$ have PUT condition and $(H^n,X^n)\Rightarrow(H,X)$, then we have  
$$(H^n, X^n,\int H^ndX^n)\Rightarrow(H,X,\int HdX).$$
\end{theorem}
\subsection{Two-level irregular discretization scheme}
Let us consider the following SDE
\begin{equation}\label{sde}
d X_t=b\left(X_t\right) d t+\sigma\left(X_t\right) d W_t, \quad X_0=x \in \mathbb{R}^d,
\end{equation}
where $W=\left(W^1 \ldots, W^q\right)$ is a $q$-dimensional Brownian motion on  $\left(\Omega, \mathcal{F},\left(\mathcal{F}_t\right)_{t \geq 0}  \mathrm{P}\right)$ with $\left(\mathcal{F}_t\right)_{t \geq 0}$ is the completion of natural filtration of  Brownian filtration. $b$ and $\sigma$ are, respectively, $\mathbb{R}^d$ and $\mathbb{R}^{d \times q}$ valued functions and satisfy the Lipschitz condition. \par
The coarse sampling is defined in (\ref{time}) and we set
\[\eta_n(t)=\tau^n_k,\quad \tau^n_k\le t<\tau^n_{k+1}.\]
Moreover, the finer sampling has the following form: fix $k$, for $0\le j\le m-1$, we let $\tau^{nm}_{k,0}=\tau^n_k$, $\tau^{nm}_{k,m}=\tau^n_{k+1}$ and 
\[\tau^{nm}_{k,j+1}=\tau^{nm}_{k,j}+\frac{1}{nm\theta(\tau^{nm}_{k,j})},\]
then we can define similarly
\[\eta_{nm}(t)=\tau^{nm}_{k,j},\quad \tau^{nm}_{k,j}\le t<\tau^{nm}_{k,j+1}.\]
We rewrite (\ref{sde}) as follows 
\begin{equation}\label{rewr}
d X_t=\varphi\left(X_t\right) d Y_t=\sum_{j=0}^q \varphi_j\left(X_t\right) d Y_t^j 
\end{equation}
where $\varphi_j$ is the $j$ the column of the matrix $\sigma$, for $1 \leq j \leq q, \varphi_0=b$ and $Y_t:=$ $\left(t, W_t^1, \ldots, W_t^q\right)^{\prime}$. We define the coarse and the finer irregular discretization scheme 
\begin{equation}
d X_t^n=\varphi\left(X_{\eta_n(t)}^n\right) d Y_t=\sum_{j=0}^q \varphi_j\left(X_{\eta_n(t)}^n\right) d Y_t^j,
\end{equation}
\begin{equation}
d X_t^{nm}=\varphi\left(X_{\eta_{nm}(t)}^{nm}\right) d Y_t=\sum_{j=0}^q \varphi_j\left(X_{\eta_{nm}(t)}^{nm}\right) d Y_t^j . 
\end{equation}
The next lemma is important to verify the conditions for the convergence in terms of quadratic variation and covariation processes and can be found in \cite{lindberg2013error}.
\begin{lemma}\label{carl}
Assume $\theta,\tau^n_k,\eta_n$ are defined above. Set 
\[\psi_n(t)=n^{\frac{p}{2}}\sum_{k=1}^{\infty}(W_t-W_{\tau^n_k})^p1_{\{\tau^n_k\le t <\tau^n_{k+1}\}},\]
and if $\theta$ is a.s. Riemann integrable over $[0,T]$ and is always positive, then for $p=1,2$,
\[\sup_{0\le t\le T}|\int_0^t \psi_n(s)ds- a_p\int_0^t\frac{1}{\theta(s)^{\frac{p}{2}}}ds|\stackrel{\mathbb{P}}{\longrightarrow} 0,\]
where $a_1$ equals $0$ and $a_2$ equals $\frac{1}{2}$.
\end{lemma}

\subsection{Multilevel Monte Carlo method}
Let $(X_t^{m^{\ell}})_{0\le t\le T}$ be the scheme defined in the previous subsection with upscript $m^{\ell}$ for $ \ell\in\{0,\cdots,L\}$, where $L=\log n/\log m$. Our task is to compute the expected value $\mathbb{E}[f(X_T)]$ for various test functions. We first regard $\mathbb{E}[f(X_T^n)]$ as an estimator and we can write $\mathbb{E}[f(X_T^n)]$ as 
\begin{equation}\label{aprox}
\mathbb{E}[f(X_T^n)]=\mathbb{E}[f(X_T^1)]+\sum_{\ell=1}^L \mathbb{E}[f(X_T^{m^{\ell}})-f(X_T^{m^{\ell-1}})].
\end{equation}
The expectation on the right-hand side of (\ref{aprox}) then can be approximated by the classical Monte Carlo method, that is
\begin{equation}\label{monte}
Q_n=\frac{1}{N_0}\sum_{k=1}^{N_0}f(X_{T,k}^1)+\sum_{\ell=1}^L\frac{1}{N_{\ell}}\sum_{k=1}^{N_{\ell}}(f(X_{T,k}^{m^{\ell}})-f(X_{T,k}^{m^{\ell-1}})),
\end{equation}
where $(X_{T,k}^{m^{\ell}},X_{T,k}^{m^{\ell-1}})_{1\le k\le N_{\ell}}$ are independent copies of $(X_{T}^{m^{\ell}},X_{T}^{m^{\ell-1}})$ and $(X_{T,k}^1)_{1\le k\le N_0}$ are independent copies of $X_T^1$. The central limit theorem for the error $Q_n-\mathbb{E}[f(X_T)]$  will be given as soon as the main theorem in the next section is proved. We should also notice that without loss of generality, we can study the scheme $( X_t^{nm})$ instead of $(X_t^{m^{\ell}})$ if we let $n={m^{\ell-1}}$.

\section{Main result}
We will prove a stable convergence theorem for the error $X_t^{nm}-X_t^n$. This type of convergence can be found in 
\cite{ben2015central} Theorem 3, where the time step is constant. Our theorem is an extension in the present setting since the irregular discretization scheme has the form of (\ref{time}). Here, we give the main result. 
\begin{theorem}\label{main}
Assume that $\theta,X^{nm},X^n$ are defined as in Subsection 2.2 and  that $b$ and $\sigma$ are $\mathcal{C}^1$ with linear growth, then we have for $m\in \mathbb{N}\backslash\{0,1\}$,
\[\sqrt{\frac{mn}{(m-1)}}(X^{nm}-X^n)\stackrel{s-\mathcal{L}}{\longrightarrow}U\quad \text{as $n\to \infty$,}\]
with $\left(U_t\right)_{0 \leq t \leq T}$ the d-dimensional process satisfying
$$
U_t=\frac{1}{\sqrt{2}} \sum_{i, j=1}^q Z_t \int_0^t \frac{H_s^{i, j}}{\sqrt{\theta(s)}} d B_s^{i j}, \quad t \in[0, T],
$$
where
 $$H_s^{i, j}=\left(Z_s\right)^{-1} \dot{\varphi}_{s, j} \bar{\varphi}_{s, i} \quad \text{with } \dot{\varphi}_{s, j}:=\nabla \varphi_j\left(X_s\right) \text{and } \bar{\varphi}_{s, i}:=\varphi_i\left(X_s\right),$$ and $\left(Z_t\right)_{0 \leq t \leq T}$ is the $\mathbb{R}^{d \times d}$ valued process solution of the linear equation
$$
Z_t=I_d+\sum_{j=0}^q \int_0^t \dot{\varphi}_{s, j} d Y_s^j Z_s, \quad t \in[0, T] .
$$
Here, $\nabla \varphi_j$ is a $d \times d$ matrix with $\left(\nabla \varphi_j\right)_{i k}$ is the partial derivative of $\varphi_{i j}$ with respect to the $k$ the coordinate, and $\left(B^{i j}\right)_{1 \leq i, j \leq q}$ is a standard $q^2$-dimensional Brownian motion independent of $W$. This process is defined on an extension $(\widetilde{\Omega}, \widetilde{\mathcal{F}},(\widetilde{\mathcal{F}}_t)_{t \geq 0}, \widetilde{\mathbb{P}})$ of the space $\left(\Omega, \mathcal{F},\left(\mathcal{F}_t\right)_{t \geq 0}, \mathbb{P}\right)$.
\end{theorem}
\begin{remark}
If we set $\theta=1/T$ in (\ref{time}), then we are in the case of the Euler scheme defined in \cite{ben2015central} with time step $T/n$. Moreover, we can recover the same limit as in Theorem 3 in \cite{ben2015central}. Hence, our result is a generalization of that theorem.
\end{remark}
To prove the theorem, we first proceed to obtain the following lemma.
\begin{lemma}\label{yux}
For $i,j\in \{1,\cdots,q\}^2$, we have
\begin{equation}
\sup_{0\le t\le T}|\sqrt{n}\int_0^t (s-\eta_n(s))dY_s^j|\stackrel{\mathbb{P}}{\longrightarrow}0,
\end{equation}

\begin{equation}
\sqrt{n}\int_0^{\cdot}(Y_s^i-Y^i_{\eta_n(s)})dY_s^j \stackrel{s-\mathcal{L}}{\longrightarrow}\frac{B^{ij}}{\sqrt{2\theta}},
\end{equation}
where $\left(B^{i j}\right)_{1 \leq i, j \leq q}$ is a standard $q^2$-dimensional Brownian motion independent of $W$. This process is defined on an extension $(\widetilde{\Omega}, \widetilde{\mathcal{F}},(\widetilde{\mathcal{F}}_t)_{t \geq 0}, \widetilde{\mathbb{P}})$ of the space $\left(\Omega, \mathcal{F},\left(\mathcal{F}_t\right)_{t \geq 0}, \mathbb{P}\right)$. Moreover, these two limiting processes satisfy PUT condition.
\end{lemma}
{\bf  Proof}. Let 
\[M^n=\sqrt{n}\int_0^t (s-\eta_n(s))dY_s^j,\]
then $M^n$ is a square-integrable martingale, and its quadratic variation process is 
\[\left\langle M^n\right\rangle_T=n\int_0^T(s-\eta_n(s))^2ds.\]
From the definition of the $\eta_n$ and the bounded process $1/\theta$, we get
\[\left\langle M^n\right\rangle_T\le\frac{TK^2}{n},\]
so we deduce that $\sup_n\mathbb{E}[\left\langle M^n   \right\rangle_T]<\infty,$ which is the definition of the PUT condition and
\[\left\langle M^n\right\rangle_T\stackrel{\mathbb{P}}{\longrightarrow} 0 \quad \text{as n $\to \infty$},\]
then we complete the first result. (9) is a direct result of Theorem 3.3 in  \cite{lindberg2013error} for the second convergence. In effect, in that paper, $ Y$ is  the solution of the SDE
\[dY_t=\alpha(Y_t)dt+\beta(Y_t)dW_t,\]
where W is a d-dimensional Brownian motion, $\alpha$ and $\beta$ are, respectively, $\mathbb{R}^d$ and $\mathbb{R}^{d \times d}$ valued functions and satisfy the global Lipschitz condition, but here  we let $d=q$, $\alpha=0$ and  $\beta=I^q$, then we deduce from the  \cite{lindberg2013error} Theorem 3.3 that 
\[Z_{i,j}^n=\sqrt{n}\int_0^t(Y_s^i-Y^i_{\eta_n(s)})dY_s^j \stackrel{s-\mathcal{L}}{\longrightarrow}\sum_{r,k=1}^q\frac{\beta_{j,k}\beta_{i,r}}{\sqrt{2\theta}}\cdot B^{r,k}=\frac{B^{i,j}}{\sqrt{2\theta}},\]
which is the second convergence and  implies the PUT condition.$\hfill\square$\\

{\bf  Proof of the Theorem \ref{main}}.

Step 1: Decomposition for the error process. 

This standard step can be found in \cite{kurtz1991wong} and in \cite{ben2015central}. Now let us consider the error process $U^{nm, n}=\left(U_t^{nm, n}\right)_{0 \leq t \leq T}$, defined by
$$
U_t^{nm, n}:=X_t^{nm}-X_t^n, \quad t \in[0, T].
$$
Through the Taylor expansion, we can easily get
$$
d U_t^{nm, n}=\sum_{j=0}^q \dot{\varphi}_{t, j}^n\left(X_{\eta_{nm}(t)}^{nm}-X_{\eta_n(t)}^n\right) d Y_t^j,
$$
where $\dot{\varphi}_{t, j}^n$ is the $d \times d$ matrix whose $i$ the row is the gradient of the real-valued function $\varphi_{i j}$ at a point between $X_{\eta_n(t)}^n$ and $X_{\eta_{nm}(t)}^{nm}$. Set
$$
G_t^{nm, n}=\int_0^t \sum_{j=0}^q \dot{\varphi}_{s, j}^n\left(X_s^n-X_{\eta_n(s)}^n\right) d Y_s^j-\int_0^t \sum_{j=0}^q \dot{\varphi}_{s, j}^n\left(X_s^{nm}-X_{\eta_{nm}(s)}^{nm}\right) d Y_s^j ,
$$
then we can rewrite $U_t^{nm, n}$ as 
$$
U_t^{nm, n}=\int_0^t \sum_{j=0}^q \dot{\varphi}_{s, j}^n U_s^{nm, n} d Y_s^j+G_t^{nm, n}.
$$
To solve $U_t^{nm,n}$, in virtue of Theorem 56  in \cite{plotter2003stochastic}, let firstly $Z_t^{nm,n}$ be the $\mathbb{R}^{d\times d}$ valued solution of 
$$
Z_t^{nm, n}=I_d+\int_0^t\left(\sum_{j=0}^q \dot{\varphi}_{s, j}^n d Y_s^j\right) Z_s^{nm, n} .
$$
Then $U_t^{nm,n}$ is given by
$$
\begin{aligned}
U_t^{nm, n}=Z_t^{nm, n} [& \int_0^t\left(Z_s^{nm, n}\right)^{-1} d G_s^{nm, n} \\
& -\int_0^t\left(Z_s^{nm, n}\right)^{-1} \sum_{j=1}^q\left(\dot{\varphi}_{s, j}^n\right)^2\left(X_s^n-X_{\eta_n(s)}^n\right) d s \\
& +\int_0^t\left(Z_s^{nm, n}\right)^{-1} \sum_{j=1}^q\left(\dot{\varphi}_{s, j}^n\right)^2\left(X_s^{nm}-X_{\eta_{nm}(s)}^{nm}\right) ds].
\end{aligned}
$$
Since the increments of the Euler scheme satisfy
$$
X_s^n-X_{\eta_n(s)}^n=\sum_{i=0}^q \bar{\varphi}_{s, i}^n\left(Y_s^i-Y_{\eta_n(s)}^i\right)
$$
and
$$
X_s^{nm}-X_{\eta_{nm}(s)}^{nm}=\sum_{i=0}^q \bar{\varphi}_{s, i}^{nm}\left(Y_s^i-Y_{\eta_{nm}(s)}^i\right),
$$
with $\bar{\varphi}_{s, i}^n=\varphi_i\left(X_{\eta_n(s)}^n\right)$ and $\bar{\varphi}_{s, i}^{nm}=\varphi_i\left(X_{\eta_{nm}(s)}^{nm}\right)$, it is easy to check that
$$
\begin{aligned}
U_t^{nm, n}= & \sum_{i, j=1}^q Z_t^{nm, n} \int_0^t H_s^{i, j, nm, n}\left(Y_s^i-Y_{\eta_m(s)}^i\right) d Y_s^j+R_{t, 1}^{nm, n}+R_{t, 2}^{nm, n} \\
& -\sum_{i, j=1}^q Z_t^{nm, n} \int_0^t \tilde{H}_s^{i, j, nm, n}\left(Y_s^i-Y_{\eta_{nm}(s)}^i\right) d Y_s^j-\widetilde{R}_{t, 1}^{nm, n}-\widetilde{R}_{t, 2}^{nm, n}
\end{aligned}
$$
with
\begin{align}
\label{R1} R_{t, 1}^{nm, n}&=\sum_{i=0}^q Z_t^{nm, n} \int_0^t K_s^{i, nm, n}\left(Y_s^i-Y_{\eta_n(s)}^i\right) d s, \\
\label{R2} R_{t, 2}^{nm, n}&=\sum_{j=1}^q Z_t^{nm, n} \int_0^t H_s^{0, j, nm, n}\left(s-\eta_n(s)\right) d Y_s^j,\\
\label{R1t}\tilde{R}_{t, 1}^{nm, n} & =\sum_{i=0}^q Z_t^{nm, n} \int_0^t \tilde{K}_s^{i, nm, n}\left(Y_s^i-Y_{\eta_{nm}(s)}^i\right) d s, \\
\label{R2t}\tilde{R}_{t, 2}^{nm, n} & =\sum_{j=1}^q Z_t^{nm, n} \int_0^t \tilde{H}_s^{0, j, nm, n}\left(s-\eta_{nm}(s)\right) d Y_s^j,
\end{align}
where, for $(i, j) \in\{0, \ldots, q\} \times\{1, \ldots, q\}$,
$$
\begin{aligned}
K_s^{i, nm, n} & =\left(Z_s^{nm, n}\right)^{-1}\left(\dot{\varphi}_{s, 0}^n \bar{\varphi}_{s, i}^n-\sum_{j=1}^q\left(\dot{\varphi}_{s, j}^n\right)^2 \bar{\varphi}_{s, i}^n\right), \\
H_s^{i, j,nm, n} & =\left(Z_s^{nm, n}\right)^{-1} \dot{\varphi}_{s, j}^n \bar{\varphi}_{s, i}^n,
\end{aligned}
$$
and
$$
\begin{aligned}
\tilde{K}_s^{i, nm, n} & =\left(Z_s^{nm, n}\right)^{-1}\left(\dot{\varphi}_{s, 0}^n \bar{\varphi}_{s, i}^{nm}-\sum_{j=1}^q\left(\dot{\varphi}_{s, j}^n\right)^2 \bar{\varphi}_{s, i}^{nm}\right), \\
\tilde{H}_s^{i, j, nm, n} & =\left(Z_s^{nm, n}\right)^{-1} \dot{\varphi}_{s, j}^n \bar{\varphi}_{s, i}^{nm} .
\end{aligned}
$$
Now, let us introduce
$$
Z_t=I_d+\int_0^t \sum_{j=0}^q\left(\dot{\varphi}_{s, j} d Y_s^j\right) Z_s \quad \text { with } \dot{\varphi}_{t, j}=\nabla \varphi_j\left(X_t\right) .
$$
Then thanks to the Theorem 48 in \cite{plotter2003stochastic}, we can get the explicit form of the $\left(Z_t\right)^{-1}$ and $\left(Z_t^{nm,n}\right)^{-1}$
$$
\left(Z_t\right)^{-1}=I_d+\int_0^t\left(Z_s\right)^{-1} \sum_{j=1}^q\left(\dot{\varphi}_{s, j}\right)^2 d s-\int_0^t\left(Z_s\right)^{-1} \sum_{j=0}^q \dot{\varphi}_{s, j} d Y_s^j .
$$
$$
\left(Z_t^{nm, n}\right)^{-1}=I_d+\int_0^t\left(Z_s^{nm, n}\right)^{-1} \sum_{j=1}^q\left(\dot{\varphi}_{s, j}^n\right)^2 d s-\int_0^t\left(Z_s^{nm, n}\right)^{-1} \sum_{j=0}^q \dot{\varphi}_{s, j}^n d Y_s^j .
$$
It is obvious that $\eta_n(t)$ and $\eta_{nm}(t)$ converge in probability to the identity, therefore we can use \cite{jacod2013} \uppercase\expandafter{\romannumeral6}.6.37 to obtain that $\dot{\varphi}_{s, j}^n$ converges in probability to the $\dot{\varphi}_{s, j}$, then taking advantage of Theorem 2.5 in Jacod and Protter  \cite{jacod98}, and according to the form of $Z_t$ and $Z_t^{nm, n}$, $\left(Z_t\right)^{-1}$ and $(Z_t^{nm, n})^{-1}$  respectively, we derive the following convergences:
\begin{equation}\label{zconv}
\sup _{0 \leq t \leq T}\left|Z_t^{nm, n}-Z_t\right| \xrightarrow{\mathbb{P}} 0 \quad \text { and } \quad \sup _{0 \leq t \leq T}\left|\left(Z_t^{nm, n}\right)^{-1}-\left(Z_t\right)^{-1}\right| \xrightarrow{\mathbb{P}} 0 .
\end{equation}
Furthermore, we can rewrite $U_t^{nm, n}$ as 
$$
U_t^{nm, n}=\sum_{i, j=1}^q Z_t^{nm, n} \int_0^t H_s^{i, j}\left(Y_{\eta_{nm}(s)}^i-Y_{\eta_n(s)}^i\right) d Y_s^j+R_t^{nm, n},
$$
with
$$
R_t^{nm, n}=R_{t, 1}^{nm, n}+R_{t, 2}^{nm, n}+R_{t, 3}^{nm, n}-\widetilde{R}_{t, 1}^{nm, n}-\widetilde{R}_{t, 2}^{nm, n}-\widetilde{R}_{t, 3}^{nm, n},
$$
where $R_{t, k}^{nm, n}$ and $\widetilde{R}_{t, k}^{nm, n}, k \in\{1,2\}$, are introduced by (\ref{R1})--(\ref{R2t}) and

\begin{align}
\label{R3}R_{t, 3}^{nm, n} & =\sum_{i, j=1}^q Z_t^{nm, n} \int_0^t\left(H_s^{i, j, nm, n}-H_s^{i, j}\right)\left(Y_s^i-Y_{\eta_n(s)}^i\right) d Y_s^j \\
\widetilde{R}_{t, 3}^{nm, n} & =\sum_{i, j=1}^q Z_t^{nm, n} \int_0^t\left(\widetilde{H}_s^{i, j,nm, n}-H_s^{i, j}\right)\left(Y_s^i-Y_{\eta_{nm}(s)}^i\right) d Y_s^j .
\end{align}

Step 2: Elimination of $R_t^{nm, n}$.

Now we start to prove that $\sup_{0\le t\le T}|\sqrt{n}R_t^{nm, n}|\stackrel{\mathbb{P}}{\longrightarrow}0$. In effect, we will separately prove that for each $k=1$, $2$, $3$, we have $\sup_{0\le t\le T}|\sqrt{n}R_{t,k}^{nm, n}|\stackrel{\mathbb{P}}{\longrightarrow}0$. Then the convergence of the $\sqrt{n}\widetilde{R}_{t, k}^{nm, n}$ can be similarly proved because $m$ can be seen as a constant. In the case of $k=1$, since we have (\ref{zconv}) and $\eta_n(t)$ and $\eta_{nm}(t)$ converge in probability to the identity, combining the form of $K_s^{i, nm, n}$ and according to the (\ref{R1}), we only need  to prove the sequence of 
\begin{equation}\label{le2}
\sqrt{n}\int_0^t\left(Y_s^i-Y_{\eta_n(s)}^i\right) d s
\end{equation}
tends to $0$ in probability for the local uniform topology and satisfies PUT condition.
However, it is easy to see 
\[\sqrt{n}(Y_s^i-Y_{\eta_n(s)}^i)=n^{\frac{1}{2}}\sum_{k=1}^{\infty}(W^i_t-W^i_{\tau^n_k})1_{\{\tau^n_k\le t <\tau^n_{k+1}\}},\]
which implies the convergence in (\ref{le2}) by Lemma \ref{carl}. Meanwhile, it is clear that
\[\sup_n\mathbb{E}[\int_0^T\sqrt{n}\left|Y_s^i-Y_{\eta_n(s)}^i\right| d s]<\infty,\]
which obtains the PUT condition. Therefore, we finish the case of $k=1$ by Theorem \ref{put}. In the case of $k=2$, because of (\ref{R2}) we should prove for $j\in \{1,\cdots,q\}$, we have
\[\sup_{0\le t\le T}|\sqrt{n}\int_0^t (s-\eta_n(s))dY_s^j|\stackrel{\mathbb{P}}{\longrightarrow}0,\]
and this sequence of processes has PUT condition. However, we have proved this conclusion in the Lemma \ref{yux}. Finally, for $k=3$, thanks to the same lemma, we have
\[\sqrt{n}\int_0^t(Y_s^i-Y^i_{\eta_n(s)})dY_s^j \stackrel{s-\mathcal{L}}{\longrightarrow}\frac{B^{ij}}{\sqrt{2\theta}},\]
they have  PUT condition. Since it is obvious that
\[\sup_{0\le t\le T}|H_t^{i, j, nm, n}-H_t^{i, j}|\stackrel{\mathbb{P}}{\longrightarrow}0,\]
combining the analysis above we deduce that $\sup_{0\le t\le T}|\sqrt{n}R_{t,3}^{nm, n}|\stackrel{\mathbb{P}}{\longrightarrow}0$ by the relation (\ref{R3}).
Step 3: Convergence.

Now, we consider the asymptotic behavior of the processes
\begin{equation}\label{process}
\sum_{i, j=1}^q Z_t^{nm, n} \int_0^t H_s^{i, j}\left(Y_{\eta_{nm}(s)}^i-Y_{\eta_n(s)}^i\right) d Y_s^j.
\end{equation}
According to (\ref{zconv}), we consider the following process:
\[M_t^{n,i,j} =\int_0^t\left(Y_{\eta_{nm}(s)}^i-Y_{\eta_n(s)}^i\right) d Y_s^j,\quad i,j\in \{1,\cdots,q\}^2.\]

By Theorem 2.1 in \cite{jacod2008continuous}, it is enough to prove:
for $(i,j,i',j')\in\{1,\cdots,q\}^4$, on the one hand, if $j\ne j'$,
\begin{equation}\label{1}
 \sqrt{n}\left\langle M^{n, i, j}, Y^{j'}\right\rangle_t \stackrel{\mathbb{P}}{\longrightarrow} 0, \text{ for all $t \in[0, T]$},
\end{equation}
\begin{equation}\label{2}
n\left\langle M^{n, i, j}, M^{n, i^{\prime}, j'}\right\rangle_t \stackrel{\mathbb{P}}{\longrightarrow} 0, \text{ for all $t \in[0, T]$},
\end{equation}
on the other hand, if $j=j'$,
\begin{equation}\label{3}
 \sqrt{n}\left\langle M^{n, i, j}, Y^{j}\right\rangle_t \stackrel{\mathbb{P}}{\longrightarrow} 0, \text{ for all $t \in[0, T]$},
\end{equation}
\begin{equation}\label{4}
n\left\langle M^{n, i, j}, M^{n, i^{\prime}, j}\right\rangle_t\stackrel{\mathbb{P}}{\longrightarrow} 0, \text{ for all $i\ne i'$ and $t \in[0, T]$},
\end{equation}
\begin{equation}\label{5}
n\left\langle M^{n, i, j}\right\rangle_t  \stackrel{\mathbb{P}}{\longrightarrow} \frac{(m-1) }{ 2m} \int_0^t\frac{1}{\theta(s)}ds,\text{ for all $i= i'$ and $t \in[0, T]$}.
\end{equation}
(\ref{1}) and (\ref{2}) can be derived easily since $Y^j$ and $Y^{j'}$ are independent of each other. For (\ref{3}), notice that
\begin{align}
 \sqrt{n}\left\langle M^{n, i, j}, Y^{j}\right\rangle_t &=\sqrt{n}\int_0^t\left(Y_{\eta_{nm}(s)}^i-Y_{\eta_n(s)}^i\right) d t\notag\\
&=\sqrt{n}\int_0^t(Y_{\eta_{nm}(s)}^i-Y^i_t)-(Y_{\eta_n(s)}^i-Y^i_t) d t,\notag
\end{align}
which implies (\ref{3}) by using theorem \ref{carl} twice.
Now we proceed to study (\ref{4}). The $L^2$ norm of the left-side in (\ref{4}) can be given by
$$
\begin{aligned}
\mathbb{E} \left[\left\langle M^{n, i, j}, M^{n, i^{\prime}, j}\right\rangle_t^2\right] & =\mathbb{E}\left[\left(\int_0^t\left(Y_{\eta_{m n}(s)}^i-Y_{\eta_m(s)}^i\right)\left(Y_{\eta_{n n}(s)}^{i^{\prime}}-Y_{\eta_n(s)}^{i^{\prime}}\right) d s\right)^2\right] \\
& =\int_0^t \int_0^t\left(\mathbb{E}\left[\left(Y_{\eta_{m n}(s)}^i-Y_{\eta_n(s)}^i\right)\left(Y_{\eta_{m n}(u)}^i-Y_{\eta_n(u)}^i\right)\right]\right)^2 d s d u \\
& =2 \mathbb{E}\left[\int_{0<s<u<t} g(s, u)^2 d s d u\right],
\end{aligned}
$$
where
$$
\begin{aligned}
g(s, u)= & \eta_{m n}(s) \wedge \eta_{m n}(u)-\eta_{m n}(s) \wedge \eta_n(u) \\
& -\eta_n(s) \wedge \eta_{m n}(u)+\eta_n(s) \wedge \eta_n(u).
\end{aligned}
$$
It is worth to notice that
$$
\eta_n(s) \leq \eta_{m n}(s) \leq s \leq \eta_n(u) \leq \eta_{m n}(u) \leq u \quad \forall s \leq \eta_n(u) .
$$
Hence, $g(s, u)=0$, for $s \leq \eta_n(u)$ and $ g(s, u)=\eta_{m n}(s)-\eta_n(s)$, for $\eta_n(u)<s<u$, then we have
$$
\begin{aligned}
\mathbb{E} \left[\left\langle M^{n, i, j}, M^{n, i^{\prime}, j}\right\rangle_t^2\right] 
& =2 \mathbb{E}\left[\int_{0<\eta_n(u)<s<u<t}\left(\eta_{m n}(s)-\eta_n(s)\right)^2 d s d u\right] \\
&\le 2\mathbb{E}\left[\frac{K^2}{n^2}\int_0^t(u-\eta_n(u))du\right]\\
&\le 2\frac{K^3}{n^3}t.
\end{aligned}
$$
which implies $n^2\mathbb{E}\left[ \left\langle M^{n, i, j}, M^{n, i^{\prime}, j}\right\rangle_t^2 \right]\le 2\frac{K^3}{n}t.$ Therefore, (\ref{4}) can be reduced by letting $n\to\infty$. It remains to prove (\ref{5}), which is the most important one, and we will use Lemma \ref{carl}. To obtain (\ref{5}), i.e. to prove
\begin{equation}\label{xin}
n\int_0^t\left(Y_{\eta_{nm}(s)}^i-Y_{\eta_n(s)}^i\right)^2 d t \stackrel{\mathbb{P}}{\longrightarrow} \frac{(m-1) }{ m} \int_0^t\frac{1}{2\theta(s)}ds,
\end{equation}
for reader's convenience, we omit the superscript $i$ and replace $Y$ by Brown motion $W$. We define
\[\psi_n(t)=n(W_t-W_{\eta_{n}(t)})^2,\]
\[\psi_{nm}(t)=nm(W_t-W_{\eta_{nm}(t)})^2.\]
By Lemma \ref{carl}, we have
\begin{equation}\label{two}
\int_0^t\psi_n(s)ds\stackrel{\mathbb{P}}{\longrightarrow} \frac{1}{2}\int_0^t\frac{1}{\theta(s)}ds,
\end{equation}
and also
\begin{equation}\label{xintwo}
 \int_0^t\psi_{nm}(s)ds\stackrel{\mathbb{P}}{\longrightarrow} \frac{1}{2}\int_0^t\frac{1}{\theta(s)}ds.
\end{equation}
Hence if we set $\psi'_n(t)=n(W_t-W_{\eta_{nm}(t)})^2$, the convergence in (\ref{xintwo}) tells us that
\[\int_0^t\psi'_n(s)ds\stackrel{\mathbb{P}}{\longrightarrow} \frac{1}{2m}\int_0^t\frac{1}{\theta(s)}ds,\]
and if it subtracts the convergence in (\ref{two}), we can get
\[\int_0^t(\psi_n(s)-\psi'_n(s))ds\stackrel{\mathbb{P}}{\longrightarrow} \frac{m-1}{2m}\int_0^t\frac{1}{\theta(s)}ds.\]
By (\ref{xin}), it is sufficient to prove
\[\int_0^t(\psi_n(s)-\psi'_n(s)-n(W_{\eta_{nm}(s)}-W_{\eta_{n}(s)})^2)ds\stackrel{\mathbb{P}}{\longrightarrow} 0.\]
But 
$$
\begin{aligned}
\psi_n(s)&=n(W_s-W_{\eta_{n}(s)})^2\\
&=n(W_s-W_{\eta_{nm}(s)})^2+n(W_{\eta_{nm}(s)}-W_{\eta_{n}(s)})^2+2n(W_s-W_{\eta_{nm}(s)})(W_{\eta_{nm}(s)}-W_{\eta_{n}(s)})\\
&=\psi'_n(t)+n(W_{\eta_{nm}(s)}-W_{\eta_{n}(s)})^2+2n(W_s-W_{\eta_{nm}(s)})(W_{\eta_{nm}(s)}-W_{\eta_{n}(s)}).
\end{aligned}
$$
Hence, we need only to show
\begin{equation}\label{xin2}
n\int_0^t(W_s-W_{\eta_{nm}(s)})(W_{\eta_{nm}(s)}-W_{\eta_{n}(s)})ds\stackrel{\mathbb{P}}{\longrightarrow} 0.
\end{equation}
We take the method similar to the proof of (\ref{4}). Firstly, The $L^2$ norm of the left side of the (\ref{xin2})  can be given by 

\begin{align}\label{fina}
 &\mathbb{E}\left[\left(\int_0^t\left(W_{s}-W_{\eta_{nm(s)}}\right)\left(W_{\eta_{m n}(s)}-W_{\eta_n(s)}\right) d s\right)^2\right] \notag\\
& =\int_0^t \int_0^t\mathbb{E}\left[\left(W_{s}-W_{\eta_{nm(s)}}\right)\left(W_{u}-W_{\eta_{nm(u)}}\right)\left(W_{\eta_{nm(s)}}-W_{\eta_{n(s)}}\right)\left(W_{\eta_{nm(u)}}-W_{\eta_{n(u)}}\right)\right] d s d u \notag\\
&=2\int_{0<s<u<t}\mathbb{E}\left[\left(W_{s}-W_{\eta_{nm(s)}}\right)\left(W_{u}-W_{\eta_{nm(u)}}\right)\left(W_{\eta_{nm(s)}}-W_{\eta_{n(s)}}\right)\left(W_{\eta_{nm(u)}}-W_{\eta_{n(u)}}\right)\right] d s d u.
\end{align}

If $s\le \eta_n(u)$, by the property of conditional expectation, we have $$\mathbb{E}\left[\left(W_{s}-W_{\eta_{nm(s)}}\right)\left(W_{u}-W_{\eta_{nm(u)}}\right)\left(W_{\eta_{nm(s)}}-W_{\eta_{n(s)}}\right)\left(W_{\eta_{nm(u)}}-W_{\eta_{n(u)}}\right)\right]=0.$$ 
If $\eta_n(u)<s\le \eta_{mn}(u)$, we can similarly derive the expectation of the above term equals to $0$. Finally, if $\eta_{mn}(u)<s< u$, we have by the independent increments property of Brownian motion,
$$
\begin{aligned}
    &\mathbb{E}\left[\left(W_{s}-W_{\eta_{nm(s)}}\right)\left(W_{u}-W_{\eta_{nm(u)}}\right)\left(W_{\eta_{nm(s)}}-W_{\eta_{n(s)}}\right)\left(W_{\eta_{nm(u)}}-W_{\eta_{n(u)}}\right)\right]\\
    &=\mathbb{E}\left[\left(W_{s}-W_{\eta_{nm(s)}}\right)\left(W_{u}-W_{\eta_{nm(u)}}\right)\left(W_{\eta_{nm(s)}}-W_{\eta_{n(s)}}\right)^2\right]\\
    &\le \frac{K}{n}\mathbb{E}\left[\left(W_{s}-W_{\eta_{nm(s)}}\right)\left(W_{u}-W_{\eta_{nm(u)}}\right)\right].
\end{aligned}
$$
Hence if we set 
$$
h(s, u)=  s \wedge u-s \wedge\eta_{nm(s)} -\eta_{nm(s)} \wedge u+\eta_{nm(s)}  \wedge\eta_{nm(u)},
$$
(\ref{fina}) is less than
\[2\frac{K}{n}\mathbb{E}\left[\int_{0<\eta_{nm(u)}<s<u<t}h(s,u)dsdu\right].\]
It is easy to verify that for $\eta_{nm(u)}<s<u$, $h(s, u)=s-\eta_{nm(u)}$.
 Hence, we have
$$
\begin{aligned}
 &\mathbb{E}\left[\left(\int_0^t\left(W_{s}-W_{\eta_{nm(s)}}\right)\left(W_{\eta_{m n}(s)}-W_{\eta_n(s)}\right) d s\right)^2\right] \\
& \le 2\frac{K}{n} \mathbb{E}\left[\int_{0<\eta_{nm(u)}<s<u<t} h(s, u)d s d u\right]\\
&=2 \frac{K}{n}\mathbb{E}\left[\int_{0<\eta_{nm(u)}<s<u<t}(s-\eta_{nm(u)})d s d u\right]\\
&\le 2\frac{K^2}{n^2m}\mathbb{E}\left[\int_0^t(u-\eta_{nm(u)})du\right]\\
&\le 2\frac{K^3}{n^3m^2}t,
\end{aligned}
$$
which yields (\ref{xin2}). Hence (\ref{5}) is proved, which completes the proof.$\hfill\square$

\section{Applications to the Monte Carlo method.}
Our task is to prove a central limit theorem for the multilevel Monte Carlo method. In the present case, we extend Theorem 4 in  \cite{ben2015central} based on the main result proved in the previous section. We take the same notation defined in subsection 2.3, and we write them again for the reader's convenience. We are expected to compute the error, which is described in the following form:
\[\varepsilon_n=\mathbb{E}[f(X_T^n)]-\mathbb{E}[f(X_T)],\]
and to approximate $=\mathbb{E}[f(X_T^n)]$, we set the estimator  defined as in (\ref{monte}):
\begin{equation}\label{rewrite}
Q_n=\frac{1}{N_0}\sum_{k=1}^{N_0}f(X_{T,k}^1)+\sum_{\ell=1}^L\frac{1}{N_{\ell}}\sum_{k=1}^{N_{\ell}}(f(X_{T,k}^{m^{\ell}})-f(X_{T,k}^{m^{\ell-1}})),
\end{equation}
It is crucial that $(X_{T,k}^{m^{\ell}},X_{T,k}^{m^{\ell-1}})_{1\le k\le N_{\ell}}$ are independent copies of $(X_{T}^{m^{\ell}},X_{T}^{m^{\ell-1}})$ and $(X_{T,k}^1)_{1\le k\le N_0}$ are independent copies of $X_T^1$. Taking advantage of the limit theorem in the above section, we can now establish a central limit theorem of Lindeberg-Feller type on the multilevel Monte Carlo method. To do so, we introduce a real sequence $\left(a_{\ell}\right)_{\ell \in \mathbb{N}}$ of positive terms such that for $p>2$
\begin{equation}\label{aell}
\lim _{L \rightarrow \infty} \sum_{\ell=1}^L a_{\ell}=\infty \quad \text { and } \quad \lim _{L \rightarrow \infty} \frac{1}{\left(\sum_{\ell=1}^L a_{\ell}\right)^{p / 2}} \sum_{\ell=1}^L a_{\ell}^{p / 2}=0.
\end{equation}
If $a_{\ell}$ satisfies these assumptions, we can use the Toeplitz lemma to get that if $\left(x_{\ell}\right)_{\ell \geq 1}$ is a sequence converging to $x \in \mathbb{R}$ as $\ell\to\infty$  then.
$$
\lim _{L \rightarrow+\infty} \frac{\sum_{\ell=1}^L a_{\ell} x_{\ell}}{\sum_{\ell=1}^L a_{\ell}}=x,
$$
we assume that the sequence of $\{N_{\ell}\}_{\ell \geq 1}$  has the  following special form, which is crucial to verify one of the conditions in Lyapunov central limit theorem:
\begin{equation}\label{nell}
N_{\ell}=\frac{n^{2 \alpha}(m-1) T}{m^{\ell} a_{\ell}} \sum_{\ell=1}^L a_{\ell}, \quad \ell \in\{0, \ldots, L\} \text { and } L=\frac{\log n}{\log m} .
\end{equation}
In the sequel, we will denote by $\widetilde{\mathbb{E}}$, respectively, $\widetilde{\text { Var }}$ the expectation, respectively, the variance defined on the probability space $(\widetilde{\Omega}, \widetilde{\mathcal{F}}, \widetilde{\mathbb{P}})$ introduced in Theorem \ref{main}.
\begin{theorem}\label{second}
Assume that $b$ and $\sigma$ are $\mathcal{C}^1$ functions satisfying the  Lipschitz condition. Let $f$ be a differentiable, real-valued, and satisfying function.
$$
\left(\mathcal{H}_f\right) \quad|f(x)-f(y)| \leq C\left(1+|x|^p+|y|^p\right)|x-y| \quad \text { for some } C, p>0 \text {. }
$$
Assume that for some $\alpha \in[1 / 2,1]$ we have
\begin{equation}\label{tiaojian}
\lim _{n \rightarrow \infty} n^\alpha \varepsilon_n=C_f(T, \alpha) .
\end{equation}
Then, for the choice of $N_{\ell}, \ell \in\{0,1, \ldots, L\}$ given by (\ref{nell}), we have
$$
n^\alpha\left(Q_n-\mathbb{E}\left(f\left(X_T\right)\right)\right) \Rightarrow \mathcal{N}\left(C_f(T, \alpha), \sigma^2\right)
$$
with $\sigma^2=\widetilde{\operatorname{Var}}\left(\nabla f\left(X_T\right) \cdot U_T\right)$ and $\mathcal{N}\left(C_f(T, \alpha), \sigma^2\right)$ denotes a normal distribution.
\end{theorem}

{\bf  Proof of the Theorem \ref{second}}. 

To simplify our notation, we give the proof for $\alpha=1$, the case $\alpha \in[1 / 2,1)$ is a straightforward deduction. Because of the form of the $Q_n$ in (\ref{rewrite}), we have
$$
Q_n-\mathbb{E}\left(f\left(X_T\right)\right)=\widehat{Q}_n^1+\widehat{Q}_n^2+\varepsilon_n,
$$
where
$$
\begin{aligned}
\widehat{Q}_n^1 & =\frac{1}{N_0} \sum_{k=1}^{N_0}\left(f\left(X_{T, k}^1\right)-\mathbb{E}\left(f\left(X_T^1\right)\right)\right), \\
\widehat{Q}_n^2 & =\sum_{\ell=1}^L \frac{1}{N_{\ell}} \sum_{k=1}^{N_{\ell}}\left(f\left(X_{T, k}^{ m^{\ell}}\right)-f\left(X_{T, k}^{ m^{\ell-1}}\right)-\mathbb{E}\left(f\left(X_T^{ m^{\ell}}\right)-f\left(X_T^{ m^{\ell-1}}\right)\right)\right) .
\end{aligned}
$$
Using assumption (\ref{tiaojian}), we obtain the term $C_f(T, \alpha)$ in the limit. Taking $N_0=\frac{n^2(m-1) T}{a_0} \sum_{\ell=1}^L a_{\ell}$, we can apply the classical central limit theorem to $\widehat{Q}_n^1$. Then we have $n \widehat{Q}_n^1 \xrightarrow{\mathbb{P}} 0$. Finally, we have only to study the convergence of $n \widehat{Q}_n^2$ and we will conclude by establishing
$$
n \widehat{Q}_n^2 \Rightarrow \mathcal{N}\left(0, \widetilde{\operatorname{Var}}\left(\nabla f\left(X_T\right) \cdot U_T\right)\right) .
$$
To get this result, we need the Lyapunov central limit theorem for a triangular array, and we set
$$
X_{n, \ell} :=\frac{n}{N_{\ell}} \sum_{k=1}^{N_{\ell}} Z_{T, k}^{m^{\ell}, m^{\ell-1}}
$$
and
$$
Z_{T, k}^{m^{\ell} m^{\ell-1}}  :=f\left(X_{T, k}^{ m^{\ell}}\right)-f\left(X_{T, k}^{ m^{\ell-1}}\right)-\mathbb{E}\left(f\left(X_{T, K}^{ m^{\ell}}\right)-f\left(X_{T, k}^{ m^{\ell-1}}\right)\right) .
$$
Hence, we have to check the following conditions:
\begin{itemize}
    \item $\displaystyle \lim _{n \rightarrow \infty} \sum_{\ell=1}^L \mathbb{E}\left(X_{n, \ell}\right)^2=\widetilde{\operatorname{Var}}\left(\nabla f\left(X_T\right) \cdot U_T\right)$.
    \item  There exists $p>2$ such that $ \displaystyle\lim _{n \rightarrow \infty} \sum_{\ell=1}^L \mathbb{E}\left|X_{n, \ell}\right|^p=0$.
\end{itemize}

For the first one, we have
\begin{align}\label{diyi}
\sum_{\ell=1}^L \mathbb{E}\left(X_{n, \ell}\right)^2 & =\sum_{\ell=1}^L \operatorname{Var}\left(X_{n, \ell}\right) \notag\\
& =\sum_{\ell=1}^L \frac{n^2}{N_{\ell}} \operatorname{Var}\left(Z_{T, 1}^{m^{\ell}, m^{\ell-1}}\right)\notag \\
& =\frac{1}{\sum_{\ell=1}^L a_{\ell}} \sum_{\ell=1}^L a_{\ell} \frac{m^{\ell}}{(m-1) T} \operatorname{Var}\left(Z_{T, 1}^{m^{\ell}, m^{\ell-1}}\right) .
\end{align}
applying the Taylor expansion theorem twice, we get
$$
\begin{aligned}
f\left(X_T^{ m^{\ell}}\right) & -f\left(X_T^{ m^{\ell-1}}\right) \\
= & \nabla f\left(X_T\right) \cdot \left(X_T^{ m^{\ell}}-X_T^{ m^{\ell-1}}\right)+\left(X_T^{ m^{\ell}}-X_T\right) \varepsilon\left(X_T, X_T^{ m^{\ell}}-X_T\right) \\
& -\left(X_T^{ m^{\ell-1}}-X_T\right) \varepsilon\left(X_T, X_T^{ m^{\ell-1}}-X_T\right) .
\end{aligned}
$$
The function $\varepsilon$ is given by the Taylor-Young expansion, so it satisfies 
$$\varepsilon\left(X_T, X_T^{ m^{\ell}}-X_T\right)\stackrel{\mathbb{P}}{\longrightarrow} 0\quad \text{ and }\quad  \varepsilon\left(X_T, X_T^{\ell, m^{\ell-1}}-X_T\right)\stackrel{\mathbb{P}}{\longrightarrow} 0$$ as $\ell\to\infty$. In effect, through the similar (but simpler) proof of Theorem \ref{main} we can get that $\sqrt{\frac{m^{\ell}}{(m-1) T}}\left(X_T^{\ell, m^{\ell}}-X_T\right)$ and $\sqrt{\frac{m^{\ell}}{(m-1) T}}\left(X_T^{\ell, m^{\ell-1}}-X_T\right)$ converge and then we deduce as $\ell\to\infty$,
$$
\begin{aligned}
\sqrt{\frac{m^{\ell}}{(m-1) T}} & \left(\left(X_T^{\ell, m^{\ell}}-X_T\right) \varepsilon\left(X_T, X_T^{\ell, m^{\ell}}-X_T\right)\right. \\
& \left.-\left(X_T^{\ell, m^{\ell-1}}-X_T\right) \varepsilon\left(X_T, X_T^{\ell, m^{\ell-1}}-X_T\right)\right) \stackrel{\mathbb{P}}{\longrightarrow}0 .
\end{aligned}
$$
So, according to Theorem \ref{main} we conclude that
$$
\sqrt{\frac{m^{\ell}}{(m-1) T}}\left(f\left(X_T^{\ell, m^{\ell}}\right)-f\left(X_T^{\ell, m^{\ell-1}}\right)\right) \Rightarrow^{\text {stably }} \nabla f\left(X_T\right) \cdot U_T \\
$$
as $ \ell \rightarrow \infty $. Moreover, if we consider relation (\ref{process}) for the decomposition of $(X_T^{\ell, m^{\ell}}-X_T^{\ell, m^{\ell-1}})$,  in virtue of BDG inequality, we can derive immediately that for $p>2$,
$$
 \sup _{\ell} \mathbb{E}\left|\sqrt{\frac{m^{\ell}}{(m-1) T}} (X_T^{\ell, m^{\ell}}-X_T^{\ell, m^{\ell-1}})\right|^p<\infty ,
$$
which implies 
\begin{equation}\label{cuo}
\sup _{\ell} \mathbb{E}\left|\sqrt{\frac{m^{\ell}}{(m-1) T}}\left(f\left(X_T^{\ell, m^{\ell}}\right)-f\left(X_T^{\ell, m^{\ell-1}}\right)\right)\right|^p<\infty
\end{equation}
Then, we deduce that
$$
\mathbb{E}\left(\sqrt{\frac{m^{\ell}}{(m-1) T}}\left(f\left(X_T^{\ell, m^{\ell}}\right)-f\left(X_T^{\ell, m^{\ell-1}}\right)\right)\right)^2 \rightarrow \widetilde{\mathbb{E}}\left(\nabla f\left(X_T\right) \cdot U_T\right)^2<\infty
$$
Consequently,
$$
\frac{m^{\ell}}{(m-1) T} \operatorname{Var}\left(Z_{T, 1}^{m^{\ell}, m^{\ell-1}}\right) \rightarrow \widetilde{\operatorname{Var}}\left(\nabla f\left(X_T\right) \cdot U_T\right)<\infty .
$$
Hence, combining this result with relation (\ref{diyi}), we obtain the first condition using Toeplitz lemma. Concerning the second one, by Burkhölder's inequality and elementary computations, we get for $p>2$
$$
\mathbb{E}\left|X_{n, \ell}\right|^p=\frac{n^p}{N_{\ell}^p} \mathbb{E}\left|\sum_{\ell=1}^{N_{\ell}} Z_{T, 1}^{m^{\ell}, m^{\ell-1}}\right|^p \leq C_p \frac{n^p}{N_{\ell}^{p / 2}} \mathbb{E}\left|Z_{T, 1}^{m^{\ell}, m^{\ell-1}}\right|^p,
$$
where $C_p$ is a numerical constant depending only on $p$. By  (\ref{cuo}) again we get that there exists a constant $K_p>0$ such that
$$
\mathbb{E}\left|Z_{T, 1}^{m^{\ell}, m^{\ell-1}}\right|^p \leq \frac{K_p}{m^{p \ell / 2}} .
$$
Therefore,
$$
\begin{aligned}
\sum_{\ell=1}^L \mathbb{E}\left|X_{n, \ell}\right|^p & \leq \tilde{C}_p \sum_{\ell=1}^L \frac{n^p}{N_{\ell}^{p / 2} m^{p \ell / 2}} \\
& \leq \frac{\tilde{C}_p}{\left(\sum_{\ell=1}^L a_{\ell}\right)^{p / 2}} \sum_{\ell=1}^L a_{\ell}^{p / 2} \underset{n \rightarrow \infty}{\longrightarrow} 0 .
\end{aligned}
$$
This completes the proof.$\hfill\square$

\section{Application to call option.}
Here, our setting is the simplest Black-Scholes model, and the market is complete, which means every option can be hedged by a pair of portfolios. We will give the limit distribution of the Black-Scholes hedging error based on the two-level sampling scheme. Assume our stock price $X^1$ is determined by a geometric Brown motion of the form
\[dX^1_t=\mu X^1_tdt+\sigma X^1_t dB_t,\]
while the bank account $X^2$ is modeled by
\[dX^2_t=rX^2_tdt.\]
Let the call option with a pay-off function $\max\{X^1_T,0\}$ for some constant $T$ and $K$ be defined as 
\begin{equation}\label{opt}
V_t:=\Phi\left(d_{+}\right) X^1_t-K e^{-(T-t)} \Phi\left(d_{-}\right),
\end{equation}
where $\Phi$ denotes the standard normal cumulative distribution function and
$$
d_{ \pm}(t)=\frac{\log (X^1_t / K)+\left(r \pm \sigma^2 / 2\right)(T-t)}{\sigma \sqrt{T-t}} .
$$
Now, if we set
$$
X_t=\binom{X^1_t}{X^2_t}.
$$
Then we have $d=2$ and $q=1$ in the Subsection 2.2. Moreover, the $b$ and $\sigma$ in the (\ref{sde}) are defined as follows:
\[b(X_t)=\binom{\mu X^1_t}{rX^2_t},\quad \text{and }\quad \sigma(X_t)=\binom{\sigma X^1_t}{0}.\]
Then we rewrite it to the same form of (\ref{rewr}) 
\[d X_t=\varphi\left(X_t\right) d Y_t=\sum_{j=0}^1 \varphi_j\left(X_t\right) d Y_t^j,\quad t\in[0,T] \]
with $\varphi_0=b(X_t)$, $\varphi_1=\sigma(X_t)$. Note that for every $t\in[0,T]$, $V_t$ in (\ref{opt}) can be seen as a function of $X_t$ having the following form
\[V_t=f(X_t)=\Phi\left(d_{+}\right) X^1_t-K e^{-T}X^2_t \Phi\left(d_{-}\right),\]
The limiting distribution of the Black-Scholes hedging error can be obtained immediately by applying the Theorem \ref{main} and Taylor expansion.
\begin{theorem}
Assume we are in the case of Theorem \ref{main}. Then if $V^{nm}_t=f(X^{nm}_t)$ and $V^n_t=f(X^n_t)$, we have 
\[\sqrt{\frac{mn}{(m-1)}}(V^{nm}_t-V^n_t)\Rightarrow \langle f',U\rangle\quad \text{for all t, as $n\to \infty$}\]
with $f'$ the derivative of f satisfying
\[f'=\binom{\Phi'\left(d_{+}\right)X^1_t+\Phi\left(d_{-}\right)-Ke^{-T}X^2_t\Phi'\left(d_{-}\right)}{Ke^{-T}\Phi\left(d_{-}\right)}\]
and $\left(U_t\right)_{0 \leq t \leq T}$ the 2-dimensional process satisfying
$$
U_t=\frac{1}{\sqrt{2}} Z_t \int_0^t \frac{H_s}{\sqrt{\theta(s)}} d W_s, \quad t \in[0, T],
$$
where $W$ is a standard Brown motion independent of $B$ and 
$$H_s=\left(Z_s\right)^{-1}\binom{\sigma^2X_t^1}{0},$$ 
and $\left(Z_t\right)_{0 \leq t \leq T}$ is the $\mathbb{R}^{2 \times 2}$ valued process solution of the linear equation
$$
Z_t=I_2+ \int_0^t \begin{pmatrix}
\mu & 0\\
 0 &r
\end{pmatrix}dsZ_s+\
\int_0^t \begin{pmatrix}
\sigma & 0\\
 0 &0
\end{pmatrix}dB_sZ_s\quad t \in[0, T] .
$$
These processes are defined on an extension $(\widetilde{\Omega}, \widetilde{\mathcal{F}},(\widetilde{\mathcal{F}}_t)_{t \geq 0}, \widetilde{\mathbb{P}})$ of the space $\left(\Omega, \mathcal{F},\left(\mathcal{F}_t\right)_{t \geq 0}, \mathbb{P}\right)$.

\end{theorem}

\textbf{Acknowledgment} Wang was partly supported by Shandong Provincial Natural Science Foundation (No. ZR2024MA082) and the National Natural Science Foundation of China (No. 12071257).

\end{document}